\documentclass{amsart}

\usepackage{amssymb}
\usepackage{graphicx}
\usepackage{epstopdf}
\theoremstyle{plain}

\newtheorem{thm}{Theorem}

\newtheorem{cor}[thm]{Corollary}
\newtheorem{lem}[thm]{Lemma}
\newtheorem{conv}[thm]{Convention}

\theoremstyle{definition}
\newtheorem{df}[thm]{Definition}
\newtheorem{rem}[thm]{Remark}

\newtheorem{conj}[thm]{Conjecture}

\author{Michael Robinson}
\address{Center for Applied Mathematics\\657 Rhodes Hall\\Cornell
  University, Ithaca, NY 14850}
\email{robinm@cam.cornell.edu}



\title[Heteroclines of a parabolic equation]{A cell complex structure for the space of heteroclines for a semilinear parabolic equation}

\begin{document}

\begin{abstract}
It is well known that for many semilinear parabolic equations there is
a global attractor which has a cell complex structure with finite
dimensional cells.  Additionally, many semilinear parabolic equations
have equilibria with finite dimensional unstable manifolds.  In this
article, these results are unified to show that for a specific
parabolic equation on an unbounded domain, the space of heteroclinic
orbits has a cell complex structure with finite dimensional cells.
The result depends crucially on the choice of spatial dimension and
the degree of the nonlinearity in the parabolic equation, and thereby
requires some delicate treatment.
\end{abstract}

\maketitle

\section{Introduction}

In this article, the space of heteroclinic orbits of
\begin{equation}
\label{limited_pde_unstable}
\frac{\partial u(t,x)}{\partial t}=\frac{\partial^2 u(t,x)}{\partial x^2} - u^2(t,x)
+ \phi(x)
\end{equation}
is shown to have the structure of a cell complex with
finite-dimensional cells, where $u \in
C^1(\mathbb{R},C^{0,\alpha}(\mathbb{R}))$, $\phi \in L^1 \cap
C^{0,\alpha}(\mathbb{R})$, and $|\phi| \to 0$ as $|x|\to \infty$.
This result is a generalization of a well-known result that the
unstable manifolds of \eqref{limited_pde_unstable} are finite
dimensional.  Indeed Theorem 5.2.1 in \cite{Henry} can easily be made
to apply with the Banach spaces we shall choose.  Theorem 5.2.1 in
\cite{Henry} shows the existence of a smooth finite dimensional
unstable manifold locally at an equilibrium.  One can then use the
iterated time-1 map of the flow for \eqref{limited_pde_unstable} to
extend this local manifold to a maximal unstable manifold.  There are
also finite Hausdorff dimensional attractors for the forward Cauchy
problem on bounded domains \cite{JRobinson_2001}.  However, we shall
exhibit a more global approach to the finite dimensionality of the
unstable manifolds.  This approach allows us to examine the finite
dimensionality of the space of heteroclinic orbits connecting a pair
of equilibria, which is a new result in the spirit of
\cite{Floer_gradient}.  The techniques used here depend rather
delicately on both the degree of the nonlinearity (which is quadratic)
and the spatial dimension (which is 1).  Both of these are important
in the standard methodology as well, as the portion of the spectrum of
the linearization in the right half-plane needs to be bounded away
from zero.  In the case of \eqref{limited_pde_unstable}, the spectrum
in the right-half plane is discrete and consists of a finite number of
points.

Of an immediate and important concern is that there may not be any
solutions to \eqref{limited_pde_unstable} which are defined in
$C^1(\mathbb{R},C^{0,\alpha}(\mathbb{R}))$.  More particularly, are
there solutions to \eqref{limited_pde_unstable} which are defined for
all time?  This question can be answered in the affirmative
\cite{RobinsonGlobal}, so this article makes the assumption that the
space of heteroclines is nonempty.  

\section{Applications}
Equation \eqref{limited_pde_unstable} is a very simple model of
combustion.  If $\phi$ is a positive constant, then the equation
supports traveling waves.  Such traveling waves can model the
propagation of a flame through a fuel source \cite{Volpert_1994}.

In addition to a model of combustion, \eqref{limited_pde_unstable} can
also be a simple model of the population of a single species, with a
spatially-varying carrying capacity, $\phi$.  Indeed, one easily finds
that under certain conditions the behavior of solutions to
\eqref{limited_pde_unstable} is reminiscent of the growth and
(admittedly tenuous) control of invasive species
\cite{Blaustein_2001}.  It is the control of invasive species that is
of most interest, and it is also what the structure of the attaching
maps of the cell complex reveals.  In the example given in Section
\ref{example_sec}, there is one more stable equilibrium, and several
other less stable ones.  The more stable equilibrium can be thought of
as the situation where an invasive species dominates.  The task, then,
is to try to perturb the system so that it no longer is attracted to
that equilibrium.  An optimal control approach is to perturb the
system so that it barely crosses the boundary of the stable manifold
of the the undesired equilibrium, and thereby the invasive species is
eventually brought under control with minimal disturbance to the rest
of the environment.

\section{Prior work}
Equations of the form \eqref{limited_pde_unstable} have been of
interest to researchers for quite some time.  Existence and uniqueness
of solutions on short time intervals (on strips
$(-t_0,t_0)\times\mathbb{R}$) can been shown using semigroup methods
and are entirely standard \cite{ZeidlerIIA}. However, there are
obstructions to the existence of eternal solutions.  Aside from the
typical loss of regularity due to solving the backwards heat equation,
there is also a blow-up phenomenon which can spoil existence in the
forward-time solution to \eqref{limited_pde_unstable}.  Blow-up
phenonmena in the forward time Cauchy problem (where one does not
consider $t<0$) have been studied by a number of authors \cite{Fujita}
\cite{Kobayashi_1977} \cite{Weissler_1981} \cite{Klainerman_1982}
\cite{Brezis_1984} \cite{Zheng_1986} \cite{Zheng_1995}.  More
recently, Zhang {\it et al.}  (\cite{Zhang_2000} \cite{Souplet_2002}
\cite{Wrkich_2007}) studied global existence for the forward Cauchy
problem for
\begin{equation*}
\frac{\partial u}{\partial t} = \Delta u + u^p - V(x) u
\end{equation*}
for positive $u,V$.  Du and Ma studied a related problem in
\cite{DuMa2001} under more restricted conditions on the coefficients
but they obtained stronger existence results.  In fact, they found
that all of the solutions which were defined for all $t>0$ tended to
equilibrium solutions.

The boundary value problem that results from taking $x
\in \Omega \subset \mathbb{R}^n$ for some bounded $\Omega$ (instead of
$x \in \mathbb{R}^n$) has also been discussed extensively in the literature
\cite{Henry} \cite{Jost_2007} \cite{BrunovskyFiedler1989}.  For the
boundary value problem, all bounded forward Cauchy problem solutions
tend to limits as $|t|\to\infty$, and these limits are equilibrium
solutions.

Almost all of the literature (including this article) describing
eternal solutions to \eqref{limited_pde_unstable} is restricted to
discussing heteroclines.  For unbounded domains and certain symmetries
of $\phi$, one can find traveling waves.  Since the
propagation of waves in nonlinear models is of great interest in
applications, there is much written on the subject.  The general idea
is that one makes a change of variables $(t,x) \mapsto \xi = x-ct$
which reduces \eqref{limited_pde_unstable} to an ordinary differential
equation.  This ordinary differential equation describes the profile
of a traveling wave.  Powerful topologically-motivated techniques,
such as the Leray-Schauder degree, can be used to prove existence of
wave solutions to \eqref{limited_pde_unstable}.  Asymptotic methods
can be used to determine the wave speed $c$, which is often of
interest in applications.  See \cite{Volpert_1994} for a very thorough
introduction to the subject of traveling waves in
\eqref{limited_pde_unstable}.

\section{The linearization and its kernel}

We begin by considering an equilibrium solution $f$ to
\eqref{limited_pde_unstable}.  As discussed in \cite{RobinsonNonauto},
this solution has asymptotic behavior which places it in $C^2 \cap L^1
\cap L^\infty(\mathbb{R})$.  We are particularly interested in
solutions which lie in the $\alpha$-limit set of $f$, those solutions
which are defined for all $t<0$ and tend to $f$.  Center attention on
this equilibrium by applying the change of variables $u(t,x) \mapsto
u(t,x) - f(x)$ to obtain
\begin{equation}
\label{pde_unstable}
\begin{cases}
\frac{\partial}{\partial t} u(t,x) = \frac{\partial^2}{\partial x^2}
u(t,x) - 2 f(x) u(t,x) - u^2(t,x)\\ u(0,x)=h(x) \in
C^2(\mathbb{R})\\
\lim_{t\to -\infty} u(t,x) = f(x)\\
 t<0,x\in \mathbb{R}.\\
\end{cases}
\end{equation}
Thus we have a final value problem for our nonlinear equation.  All
solutions to \eqref{pde_unstable} will tend to zero as $t\to -\infty$
uniformly by Lemma 6 of \cite{RobinsonClassify}.  Of course,
\eqref{pde_unstable} is ill-posed.  We show that there is only a
finite dimensional manifold of choices of $h$ for which a solution
exists.

\subsection{Backward time decay}

The decay of solutions to zero is a crucial part of the analysis, as
it provides the ability to perform Laplace transforms.  In the forward
time direction, one obtains upper bounds for solutions by way of
maximum principles, and lower bounds for the upper bounds by way of
Harnack estimates.  In the backward time direction, these tools
reverse roles.  Harnack estimates provide upper bounds, while the maximum
principle provides lower bounds for the upper bound.  In the proof of
Lemma 6 of \cite{RobinsonClassify}, the latter was used to some
advantage.  In this section, we briefly apply a standard Harnack
estimate to obtain an exponentially decaying upper bound.

Harnack estimates for a very general class of parabolic equations are
discussed in \cite{Kurihara_1967} and \cite{Aronson_1967}.  In those articles, the authors
examine positive solutions to
\begin{equation*}
\text{div }{\bf A}(x,t,u,\nabla u) - \frac{\partial u}{\partial t} = B(x,t,u,\nabla u),
\end{equation*}
where $x\in\mathbb{R}^n$, and ${\bf
  A}:\mathbb{R}^{2n+2}\to\mathbb{R}^n$ and $B:\mathbb{R}^{2n+2}\to
  \mathbb{R}$ satisfy
\begin{eqnarray*}
|{\bf A}(x,t,u,p)| &\le& a|p| + c|u| + e\\
|B(x,t,u,p)| &\le& b|p| + d|u| + f\\
p \cdot {\bf A}(x,t,u,p) &\ge& \frac{1}{a} |p|^2 - d|u|^2 - g,\\
\end{eqnarray*}
for some $a>0$ and $b,...g$ are measurable functions.  For a solution
$u$ defined on a rectangle $R$, the authors define a pair of
congruent, disjoint closed rectangles $R^+,R^- \subset R$ with $R^-$
being a backward time translation of $R^+$.  The main result is the
Harnack inequality
\begin{equation}
\label{harnack_eqn}
\max_{R^-} u \le \gamma \left( \min_{R^+}u + L \right ),
\end{equation}
where $\gamma>0$ depends only on geometry and $a$ (but not $b,...g$)
and $L$ is a linear combination of $e,f,g$ whose coefficients depend
on geometry.  

In the case of \eqref{pde_unstable}, or indeed of the analogous
equation with higher degree terms, we have that \eqref{harnack_eqn}
will apply with $L=0$.  Notice that the conditions on $A,B$ are
satisfied because any solution to \eqref{pde_unstable} is
automatically a finite energy solution, and therefore is bounded and
has bounded first derivatives.  The only difficulty is that
\eqref{harnack_eqn} applies for {\it positive} solutions, while
\eqref{pde_unstable} may have solutions with negative portions.  However, one can
pose the problem for the (weak) solution of
\begin{eqnarray*}
\frac{\partial |u|}{\partial t} &=& \text{sgn } (u) \left(\Delta u - u^2 - 2 f u\right)\\
&=&\Delta |u| - u|u| - 2 f |u|\\
&\ge& \Delta |u| - |u|^2 - 2 |f| |u|\\
\end{eqnarray*}
for which we only get positive solutions.  By iterating
\eqref{harnack_eqn} we have that solutions to \eqref{pde_unstable} decay
exponentially as $t\to -\infty$.  

\subsection{Topological considerations}

\begin{df}
Let $Y_a(X)$ be the subspace of $C^1(X,C^{0,\alpha}(\mathbb{R}))$ which
consists of functions which decay exponentially to zero like $e^{at}$,
where $0<\alpha \le 1$.  We define the weighted norm
\begin{equation*}
\|u\|_{Y_a} = \left\|e^{-at}\|u(t)\|_{C^{0,\alpha}(\mathbb{R})}\right\|_{C^1} 
\end{equation*}
 and the space
\begin{equation*}
Y_a(X)=\left\{u=u(t,x)\in C^1(X,C^{0,\alpha}(\mathbb{R}))|
\|u\|_{Y_a} < \infty \right\}.
\end{equation*}
In a similar way, we can define the weighted Banach space $Z_a(X)$ as a
subspace of $C^0(X,C^{0,\alpha}(\mathbb{R}))$.  It is quite important
that $Y_a$ and $Z_a$ are Banach algebras under pointwise
multiplication.
\end{df}

In light of the previous section, solutions to \eqref{pde_unstable}
are zeros of the densely defined nonlinear operator $N:Y_a((-\infty,0]) \to Z_a((-\infty,0])$ given by
\begin{equation}
\label{N_def}
N(u)=\frac{\partial u}{\partial t} - \frac{\partial^2 u}{\partial x^2}
+ u^2 + 2fu.
\end{equation}
About the zero function, the linearization of $N$ is the densely
defined linear map $L:Y_a((-\infty,0]) \to Z_a((-\infty,0])$ given by
\begin{equation}
\label{L_def}
L = \frac{\partial }{\partial t} - \frac{\partial^2 }{\partial x^2}
+ 2 f = \frac{\partial }{\partial t} - H,
\end{equation}
where we define $H=\frac{\partial^2 }{\partial x^2} - 2 f$.  Also note
that $L$ is the Frech\'et derivative of $N$, which follows from the
fact that $Y_a$ and $Z_a$ are Banach algebras.

\begin{rem}
We are using $C^{0,\alpha}(\mathbb{R})$ instead of $C^0(\mathbb{R})$
to ensure that $N$ and $L$ be densely defined.  We could use space of
continous functions which decay to zero, or the space of uniformly
continous functions equally well.
\end{rem}

\begin{conv}
We shall conventionally take $a>0$ to be smaller than the smallest
eigenvalue of $H$.
\end{conv}

We show two things: that the kernel of $L$ is finite
dimensional, and that $L$ is surjective.  These two facts enable us to
use the implicit function theorem to conclude that the space of
solutions comprising the $\alpha$-limit set of an equilibrium is a
finite dimensional submanifold of $Y_a((-\infty,0])$.

\subsection{Dimension of the kernel}

\begin{lem}
\label{eq_findim}
If $f$ is an equilibrium solution, then the operator
$L:Y_a((-\infty,0])\to
Z_a((-\infty,0])$ in \eqref{L_def} has a finite
dimensional kernel.
\begin{proof}
Notice that the operator $L$ is separable, so we try the usual
separation $h(t,x)=T(t)X(x)$.  Substituting into \eqref{L_def} gives
\begin{eqnarray*}
0&=&Lh=\left(\frac{\partial}{\partial t} - \frac{\partial^2}{\partial
  x^2} + 2 f \right)h\\
&=&T'X+T\left(-\frac{\partial^2}{\partial x^2} + 2f\right)X\\
\frac{T'}{T}&=&\frac{\left(\frac{\partial^2}{\partial x^2} - 2f\right)X}{X}=\lambda\\
\end{eqnarray*}
for some $\lambda\in\mathbb{C}$.  The separated equation for $T$
yields $T=C_x e^{\lambda t}$.  Since we are looking for the kernel of
$L$ in $Y_a \subset L^\infty(\mathbb{R}^2)$, we must conclude that
$\lambda$ must have nonnegative real part.  On the other hand, the
spectrum of $H=\left(\frac{\partial^2}{\partial x^2} - 2f\right)$ is
strictly real, so $\lambda \ge 0$.  Indeed, there are finitely many
positive possibilities for $\lambda$ each with finite-dimensional
eigenspace.  This is a standard fact about the Schr\"{o}dinger
operator $H$ since $f$ is an equilibrium.  Thus $L$ has a finite
dimensional kernel.
\end{proof}
\end{lem}

\subsection{Surjectivity of the linearization}

In order to show the surjectivity of $L$, we will construct a map
$\Gamma:Z_a((-\infty,0]) \to
  Y_a((-\infty,0])$ for which $L \circ \Gamma =
    \text{id}_{Z_a}$.  That is, we construct a right-inverse to $L$,
    noting of course that $L$ is typically not injective.
We shall derive a formula for $\Gamma$ using the Laplace
transform $v \mapsto \overline{v}$
\begin{equation*}
\overline{v}(s,x) = \int_{-\infty}^0 e^{st} v(t,x) dt,
\end{equation*}
where $\Re(s)>-a$ and $v\in Z_a((-\infty,0])$.

Since Lemma \ref{eq_findim} essentially solves \eqref{pde_unstable},
we will be solving the inhomogeneous problem with zero final condition
\begin{equation}
\label{eq_findim_zero}
\begin{cases}
\frac{\partial v(t,x)}{\partial t}-\frac{\partial^2 v(t,x)}{\partial x^2} + 2
f(x) v(t,x) = -w(t,x) \in Z_a((-\infty,0])\\
v(0,x)=0\\
\end{cases}
\end{equation}
for $t<0$.  The Laplace transform of this problem is 
\begin{eqnarray*}
s\overline{v}(s,x)+\frac{\partial^2 \overline{v}(s,x)}{\partial x^2} - 2 f(x)
\overline{v}(s,x)&=&\overline{w}(s,x)\\
(H+s)\overline{v}(s,x)&=&\overline{w}(s,x).
\end{eqnarray*}
Choose a vertical contour $C$ with $0>\Re(s) > -a$, so that the
Laplace transforms are well-defined, and that the contour remains
entirely in the resolvent set of $-H$.  Then we can invert to obtain
\begin{equation*}
\overline{v}(s,x)=(H+s)^{-1}\overline{w}(s,x).
\end{equation*}
Using the inversion formula for the Laplace transform yields
\begin{eqnarray*}
v(t,x)&=& \frac{1}{2\pi i} \int_C e^{-st} (H+s)^{-1} \overline{w}(s,x)
ds\\
&=&\frac{1}{2\pi i} \int_C e^{-st} (H+s)^{-1} \int_t^0
e^{s\tau}w(\tau,x)d\tau\,ds\\
&=&\int_t^0 \left(\frac{1}{2\pi i} \int_C e^{s(\tau-t)} (H+s)^{-1} 
ds \right)w(\tau,x)d\tau.\\
\end{eqnarray*}

\begin{figure}
\begin{center}
\includegraphics[height=2.5in]{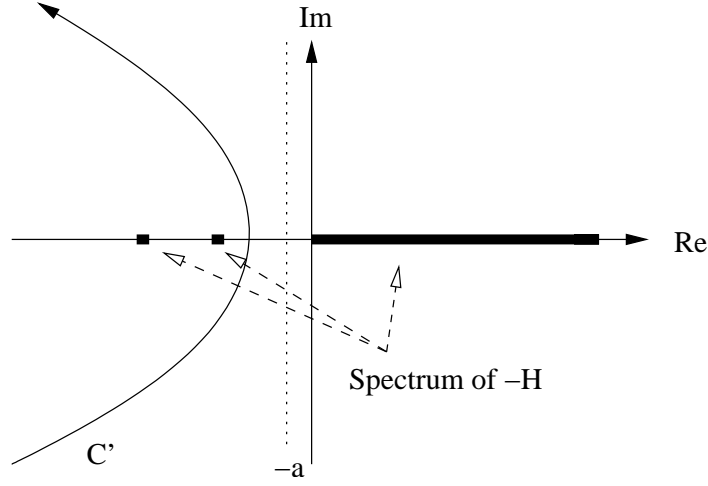}
\end{center}
\caption{Definition of the contour $C'$}
\label{Cp_contour_fig}
\end{figure}

We can obtain operator convergence of the operator-valued integral
in parentheses if we deflect the contour $C$.  Choose instead the
portion $C'$ of the hyperbola (See Figure \ref{Cp_contour_fig}) 
\begin{equation}
\left(\Re(s)\right)^2-\left(\Im(s)\right)^2 = \frac{1}{4}(\lambda-a)^2
\end{equation}
(where $\lambda$ is the smallest magnitude eigenvalue of $-H$) which
lies in the left half-plane as our new contour.  Then, since $-H:C^{0,\alpha}
\to C^{0,\alpha}$ is sectorial about $(\lambda-a)/2$, Theorem 1.3.4 in \cite{Henry} implies that
the integral
\begin{equation*}
\left(\frac{1}{2\pi i} \int_{C'} e^{s(\tau-t)} (H+s)^{-1} 
ds \right)
\end{equation*}
defines an operator-valued semigroup $e^{-H(\tau-t)}$, so the formula
  for $\Gamma$ is given by
\begin{equation}
\Gamma(w)(t,x)= \int_t^0 e^{-H(\tau-t)} w(\tau,x)d\tau.
\end{equation}
It remains to show that the image of $\Gamma$ is in fact $Y_a$,
as it is easy to see that its image is in $L^\infty$.  That the image
is as advertised is not immediately obvious because the contour
deflection $C \to C'$ changes the domain of the Laplace transform.  In
particular, the derivation given above is no longer valid with the
new contour.

Therefore, we must estimate $\|v\|_{Z_a}$ (recall that $\lambda$ is the smallest magnitude eigenvalue of $-H$)

\begin{eqnarray*}
\|e^{-at} v(t,x)\|_{C^0} &=& \left \| \frac{1}{2\pi i}\int_{C'} (s+H)^{-1}
 \int_t^0 e^{-(s+a)(t-\tau)} e^{a\tau}w(\tau,x) d\tau\,ds \right \|_{C^0}\\
&\le&\frac{1}{2\pi}\int_{C'} \frac{K_1}{|s-\lambda|}e^{-\Re(s+a)t}
 \int_t^0 e^{\Re(s+a)\tau} \|w\|_{Z_a} d\tau\,ds\\
&\le&\frac{K_1\|w\|_{Z_a}}{2\pi}\int_{C'} \frac{1}{|s-\lambda|}e^{-\Re(s+a)t}
 \frac{1}{\Re(s+a)}\left(1-e^{\Re(s+a)t}\right)ds\\
&\le&\frac{K_1\|w\|_{Z_a}}{\pi}\int_{C'}
 \frac{ds}{|s-\lambda||\Re(s+a)|}\\
&\le&K_2 \|w\|_{Z_a},\\
\end{eqnarray*}
where $0<K_1,K_2<\infty$ are independent of $t$ and $w$.  We have made
use of the usual estimate of the norm
of $(H+s)^{-1}:C^{0,\alpha} \to C^{0,\alpha}$ when $s$ is in the resolvent set of
$-H$.  In particular, note that the choice of $C'$ being to the left
of $-a$ is crucial to the convergence of the integrals.  Thus the
image of $\Gamma$ lies in $Z_a$.  The backward-time decay of
$\frac{\partial v}{\partial t}$ is immediate from the Harnack
inequality, so in fact the image of $\Gamma$ lies in $Y_a$.

\begin{thm}
\label{unstable_thm}
The linear map $L:Y_a((-\infty,0])\to
  Z_a((-\infty,0])$ is surjective and has a
  finite dimensional kernel.  Therefore the set $N^{-1}(0)$ is a
  finite dimensional manifold, which is the unstable manifold of the
  equilibrium $f$.  The dimension of $N^{-1}(0)$ is precisely the
  dimension of the positive eigenspace of $H$.
\begin{proof}
The only thing which remains to be shown is that the domain $Y_a$
splits into a pair of closed complementary subspaces: the kernel of
$L$ and its complement.  That its complement is closed follows
immediately from a standard application of the Hahn-Banach theorem.
(Extend $\text{id}_{\text{ker }L}$ to all of $Y_a$.)
\end{proof}
\end{thm}

Combining the fact that an equilibrium solution can have an
empty unstable manifold (a numerical computation of the dimension of
the eigenspaces of $L$ can be found in \cite{RobinsonNonauto}) and
is yet unstable, we have proven the following result.

\begin{thm}
All equilbrium solutions to \eqref{limited_pde_unstable} are degenerate
critical points in the sense of Morse.
\end{thm}

\section{Linearization about heteroclinic orbits}

We can extend the technique of the previous section to the
linearization about a heteroclinic orbit.  The resulting
generalization of Theorem \ref{unstable_thm} is that the connecting
manifolds of \eqref{limited_pde_unstable} are all finite dimensional.

Suppose that $u$ is a heteroclinic orbit of
\eqref{limited_pde_unstable}.  Let $f_-,f_+$ be the equilibrium
solutions of \eqref{limited_pde_unstable} to which $u$ converges as
$t\to -\infty$ and $t\to +\infty$ respectively.

Suppose that $\lambda_0:\mathbb{R} \to (0,\infty)$ is the smallest
positive eigenvalue of $H(t)$.  It is easy to see that $\lambda_0$ is
piecewise $C^1$, for instance, see Proposition I.7.2 in
\cite{Kielhoefer_2004}.  The fact that the the spectrum of $H$ lies
entirely to the left of $\max\{2 \|f_+\|_\infty, 2\|f_-\|_\infty\}$ ensures
that $\lambda_0$ is a bounded function.  We will define a pair of
bounded, piecewise $C^1$ functions $\lambda_1$ and $\lambda_2$ which
will aid us in defining a two more pairs of function spaces.  Let
$\lambda_1:\mathbb{R} \to (0,\infty)$ be a bounded, piecewise $C^1$
function with bounded derivative which has the following properties:
\begin{itemize}
\item $\lambda_1(t)$ is never an eigenvalue of $H(t)$, 
\item $\lim_{t\to\infty}\frac{\lambda_1(t)}{\lambda_0(t)}<1$,
\item $\lim_{t\to -\infty}\frac{\lambda_1(t)}{\lambda_0(t)}<1$, and
\item since $u\to
f_\pm$ uniformly, for a sufficiently large $R>0$, $\lambda_1$ can be
chosen so that there are no jumps on its restriction to $\mathbb{R} -
[-R,R]$.
\end{itemize}

Defining $\lambda_2$ is a somewhat more delicate problem.  We would
like to exclude the solutions which lie in the unstable manifold of
$f_+$, since they cannot lie in the space of heteroclines from $f_-
\to f_+$.  We do this by separating the eigenvalues corresponding to
the intersection of the unstable manifolds of $f_-$ and $f_+$ from those
which lie in the stable manifold of $f_+$.  However, there is an
obstruction to this technique.  In particular, the eigenvalues of
$H(t)=\frac{\partial^2}{\partial x^2} - 2 u(t)$ vary with time, and
can bifurcate.  To avoid this issue, we need some kind of regularity
for the eigenvalues to prevent them from bifurcating.  We follow Floer
\cite{Floer_relative} in the following way:

\begin{conj}
\label{baire_conj}
There is a generic subset (a Baire subset) of choices for $\phi$ in \eqref{limited_pde_unstable} so that if $u$ is a heteroclinic
orbit, all of the eigenvalues of $H(t)$ are simple.  
\end{conj}

Numerical evidence, as exhibited in \cite{RobinsonNonauto} and
Section \ref{example_sec} suggests that the above Conjecture is true.  When we
assume that all of the eigenvalues of $H(t)$ are simple, and therefore
do not undergo any bifurcations other than passing through zero, we
shall say $u$ is a heterocline contained in $U_{reg}$.

\begin{figure}
\begin{center}
\includegraphics[height=2.5in]{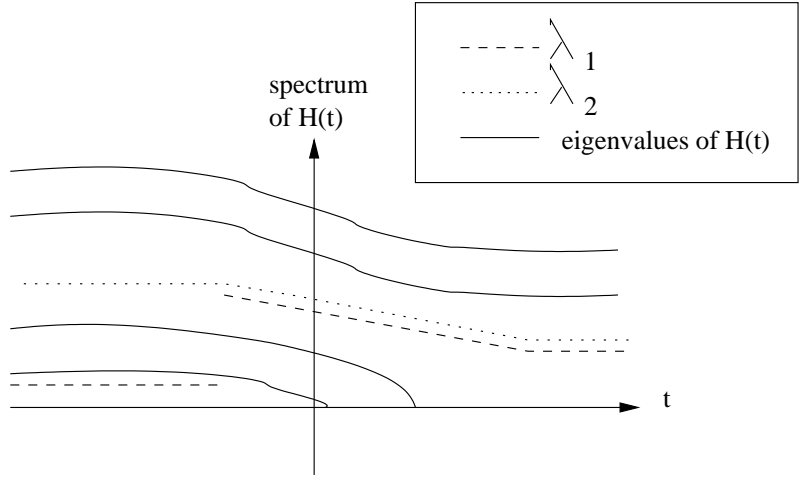}
\end{center}
\caption{Definition of $\lambda_1$ and $\lambda_2$}
\label{lambda_defs_fig}
\end{figure}

Let $\lambda_2$ be in $C^1(\mathbb{R})$ such that
\begin{itemize}
\item $\lambda_2=\lambda_1$ on $[R,\infty)$, and
\item $\lambda_2(t)$ is not an eigenvalue of
$H(t)$ for any $t$.
\end{itemize}
We can do this when $u \in U_{reg}$.  See Figure
\ref{lambda_defs_fig}.

\begin{df}
Define the Banach algebra $Y_{\lambda_i}(X)$ (for $i=1,2$) to
be the set of $u$ in $C^1(X,C^{0,\alpha}(\mathbb{R}))$ such that the
norm
\begin{equation*}
\left\|e^{-\int_0^t \lambda_i(\tau) d\tau} \|u(t)\|_{C^{0,\alpha}}
\right\|_{C^1}<\infty,
\end{equation*}
where $X$ is an interval containing zero.  Likewise, we can define the
spaces $Z_{\lambda_i}(X) \subset
C^0(X,C^{0,\alpha}(\mathbb{R}))$ in a similar way.  That these are
Banach spaces follows from the boundedness of the $\lambda_i$.  It is
also elementary to see that these are Banach algebras.
\end{df}

We then consider $N_i,L_i$ as $Y_{\lambda_i}(\mathbb{R})\to
Z_{\lambda_i}(\mathbb{R})$, where $L_i$ is the linearization of
$N_i$ about $u$ for $i=1,2$.  (Again, since $Y_{\lambda_i}$ and
$Z_{\lambda_i}$ are Banach algebras, $L_i$ is the Frech\'et
derivative of $N_i$.)  For a $i\in \{1,2\}$, consider the
restriction $L^-_i$ of $L_i$ to a map
$Y_{\lambda_i}((-\infty,0])\to
Z_{\lambda_i}((-\infty,0])$.  We rewrite
\begin{equation}
\label{lminus}
L^-_i=\left(\frac{\partial}{\partial t} - \frac{\partial^2}{\partial x^2}+2f_- \right)+(2f_- - 2u).
\end{equation}
Likewise, we can define $L^+_i:Y_{\lambda_i}([0,\infty))\to Z_{\lambda_i}([0,\infty))$.

We define the positive eigenspaces $V^+$ for the equilibria as well
\begin{equation}
V^+(f_\pm)=\text{span }\left \{v\in C^{0,\alpha}(\mathbb{R})| \text{ there
  is a }\lambda>0 \text{ with } \left(\frac{\partial^2
  }{\partial x^2} -2f_\pm \right)v = \lambda v \right \}.
\end{equation}
Note in particular that $\text{dim }V^+(f_\pm)<\infty$.

\begin{lem}
\label{findim_reqs}
If $u\in U_{reg}$ is a heterocline that converges to $f_\pm$ as
$t\to\pm\infty$, then the operator $L_i$ has a finite
dimensional kernel for $i\in\{1,2\}$, and in particular
\begin{equation*}
\lim_{t\to -\infty}\text{dim }V^+(u(t)) - \lim_{t\to +\infty}\text{dim }V^+(u(t))  \le \text{dim ker }L_i
\le\text{dim ker }L^-_i <\infty.
\end{equation*}
(The condition $u\in U_{reg}$ is only necessary for the $i=2$ case.)
\begin{proof}
Notice that the first term of \eqref{lminus} has finite dimensional
kernel by Lemma \ref{eq_findim} and closed image by Theorem \ref{unstable_thm}.
The second term of \eqref{lminus} is a compact operator since $u\to f_-$ uniformly.
Thus $L^-_i$ has a finite dimensional kernel.  Let $\text{span}\{v_m\}_{m=1}^{M}=\ker L^-_i$
and consider the set of Cauchy problems 
\begin{equation}
\label{cpset}
\begin{cases}
\frac{\partial h}{\partial t}=\frac{\partial^2 h}{\partial x^2}-2u h \text{ for }t>0\\
h(0,x)=v_m(0,x).
\end{cases}
\end{equation}
Standard parabolic theory gives uniqueness of solutions to \eqref{cpset}, and that
a solution $h$ lies in the kernel of $L^+_i$, the restriction of $L_i$ to $[0,\infty)\times\mathbb{R}$.
Therefore $\text{dim ker }L_i \le \text{dim ker }L^-_i<\infty$.  

For the other inequality, modify $u$ outside of
$[-R,R]\times\mathbb{R}$ to get a $\bar{u}$ so that the linearization $\overline{L_i}$
of $N$ about $\bar{u}$ satisfies
\begin{itemize}
\item $\text{ker }\overline{L_i}$ is isomorphic to $\text{ker }L_i$ as vector
  spaces,
\item $\bar{u}|_{(-\infty,-R)\times\mathbb{R}} = f_-$, and
\item $\bar{u}|_{(R,\infty)\times\mathbb{R}} = f_+$.
\end{itemize}
We can do this for a sufficiently large $R$, since $u$ tends uniformly
to equilibria.  Then the flow of
\begin{equation*}
\frac{\partial h}{\partial t}=\frac{\partial^2 h}{\partial x^2} +
2\bar{u} h
\end{equation*}
defines an injective linear map from the timeslice at $-R$ to the
timeslice at $R$.  (That is, it gives an injective map from
$C^{0,\alpha}(\mathbb{R})$ to itself -- injectivity being an expression
of the uniqueness of solutions.)  Each element $v$ of the kernel
of $\overline{L_i}$ evidently must have $v(-R)\in V^+(f_-)$ and
$v(R) \notin V^+(f_+)$.  Therefore, the injectivity ensures that the
intersection of the image under the flow of $V^+(f_-)$ with the
complement of $V^-(f_+)$ has at least dimension $\text{dim }V^+(f_-) -
\text{dim }V^+(f_+)$. 
\end{proof}
\end{lem}

\begin{rem}
Multiplication by $u$, $C^1(\mathbb{R}^2,C^{0,\alpha}(\mathbb{R}))\to C^0(\mathbb{R}^2)$ is not a compact operator, in particular note that $\text{dim ker }L^+_i=\infty$.
\end{rem}

\begin{thm}
\label{finite_dim_connecting_thm}
Let $u$ be a heterocline of \eqref{limited_pde_unstable}
which connects equilibria $f_\pm$.  There exists a union $\bigcup M_u$ of finite dimensional
submanifolds $M_u$ of $C^1(\mathbb{R},C^{0,\alpha}(\mathbb{R}))$ which
\begin{itemize}
\item contains $u$ and
\item consists of heteroclines connecting $f_-$ to $f_+$.
\end{itemize}
If $u\in U_{reg}$, then $M_u$ has dimension $\lim_{t\to -\infty}\text{dim }V^+(u(t)) - \lim_{t\to
  \infty}\text{dim }V^+(u(t))$, and this is maximal among such
  submanifolds $M_u$.
\begin{proof}
Observe that $L_1$ is surjective, since it is easy to show that the
formula
\begin{equation*}
\Gamma_1(w)(t)=\int_t^0 e^{-\int_0^{T-t} H(\tau) d\tau} w(T,x) dT
\end{equation*}
is a well defined right inverse of $L_1$.  This involves showing that 
\begin{equation*}
e^{-\int_0^{t} H(\tau) d\tau}=\frac{1}{2\pi i} \int_{C(t)} e^{st}
(H(t)+s)^{-1} ds
\end{equation*}
converges, where we note that the contour changes with time.  As it
happens, the computation in \cite{Henry} goes through with the only
change that at $t=0$, we deflect the contour to the right, rather than
the left (as in Figure \ref{Cp_contour_fig}).  Since Lemma
\ref{findim_reqs} shows that $L_1$ has finite dimensional kernel, then
it follows that $M_u=N_1^{-1}(0)$ is a union of finite dimensional
manifolds, with a finite maximal dimension.  It is obvious that $M_u$
consists entirely of heteroclinic orbits and contains $u$.  

It remains to show that the dimension of $M_u$ is as advertised and
maximal.  Observe that $L_2$ is a compact perturbation of an operator
$L'_2:Y_{\lambda_2}(\mathbb{R})\to
Z_{\lambda_2}(\mathbb{R})$ which is time-translation invariant.
This follows from the precise choice of $\lambda_2$ being continous
and not intersecting the eigenvalues of $H$.  $L_2$ and $L'_2$ are
both surjective by exactly the same reasoning as for $L_1$.  $L'_2$ is
injective by using separation of variables as in Lemma \ref{eq_findim}
(noting that all nontrivial solutions blow up in the
$Y_{\lambda_2}$ norm).  Therefore the Fredholm index of $L'_2$,
hence $L_2$ is zero.  However, this implies that $L_2$ is injective.

Since $L_2$ is bijective, any solution to $L_2 u = 0$ which decays
faster than $e^{\int \lambda_2(t) dt}$ as $t\to -\infty$ ends up
growing faster than $e^{\int \lambda_2(t) dt}$ as $t\to +\infty$, and
in particular does not tend to zero.  As a result, such a solution
cannot be in $\text{ker }L_1$.  This implies that $\text{dim ker }L_1
\le \lim_{t\to -\infty}\text{dim }V^+(u(t)) - \lim_{t\to
  \infty}\text{dim }V^+(u(t))$, which with the estimate in Lemma
\ref{findim_reqs} completes the proof.
\end{proof}
\end{thm}

\begin{rem}
Even if $u \notin U_{reg}$ (when there exist nonsimple eigenvalues of
$H(t)$), the function $\lambda_1$ can still be constructed.  As a
result, we {\it always} get that the connecting manifold $M_u$ is
finite-dimensional.
\end{rem}

\begin{cor}
\label{cell_complex_cor}
The space of heteroclinic orbits has the structure of a cell complex
with finite dimensional cells.  This cell complex structure is
evidently finite dimensional if there exist only finitely many
equilibria for \eqref{limited_pde_unstable}.
\end{cor}

\section{An extended example}
\label{example_sec}
\begin{figure}
\includegraphics[height=2in]{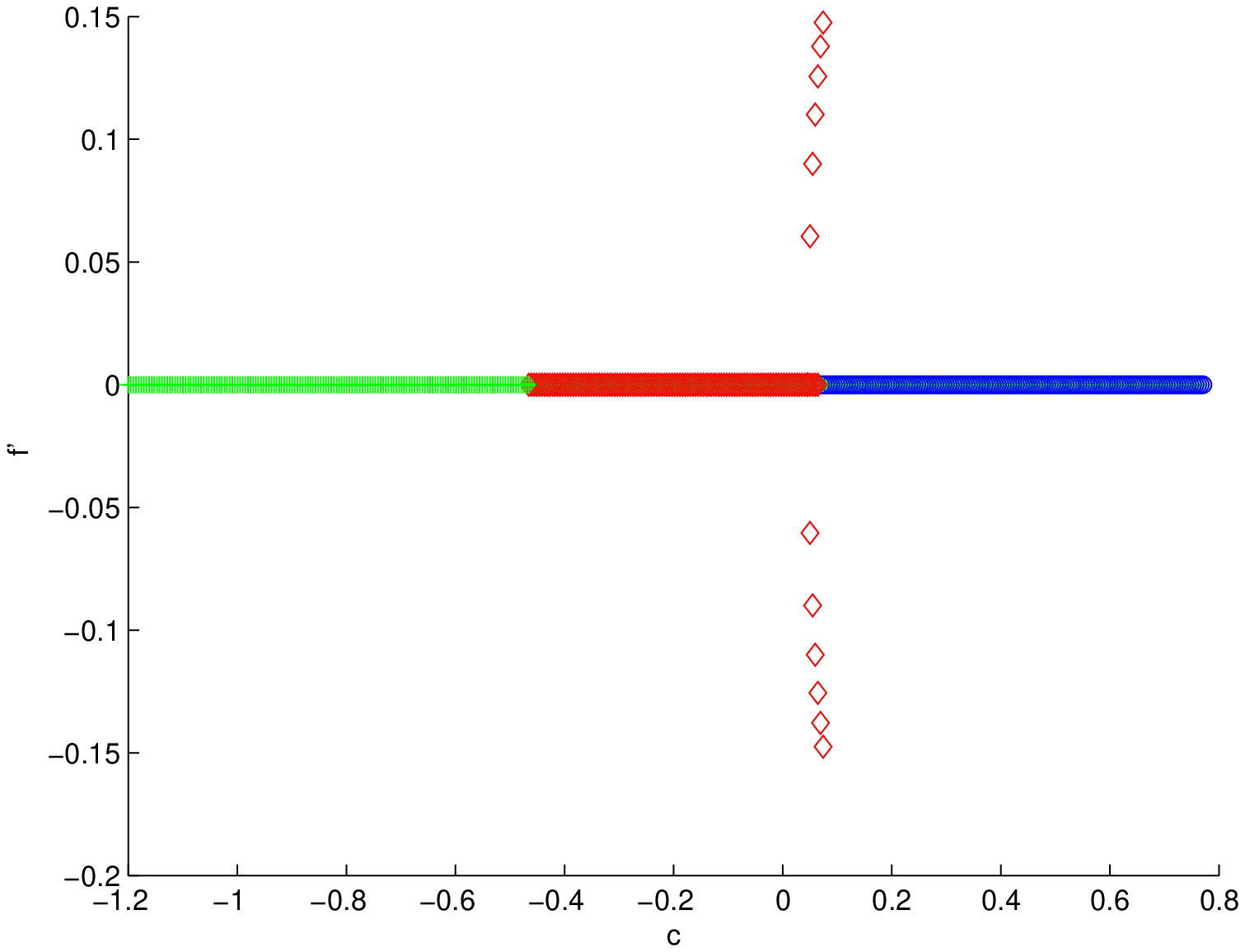}
\includegraphics[height=2in]{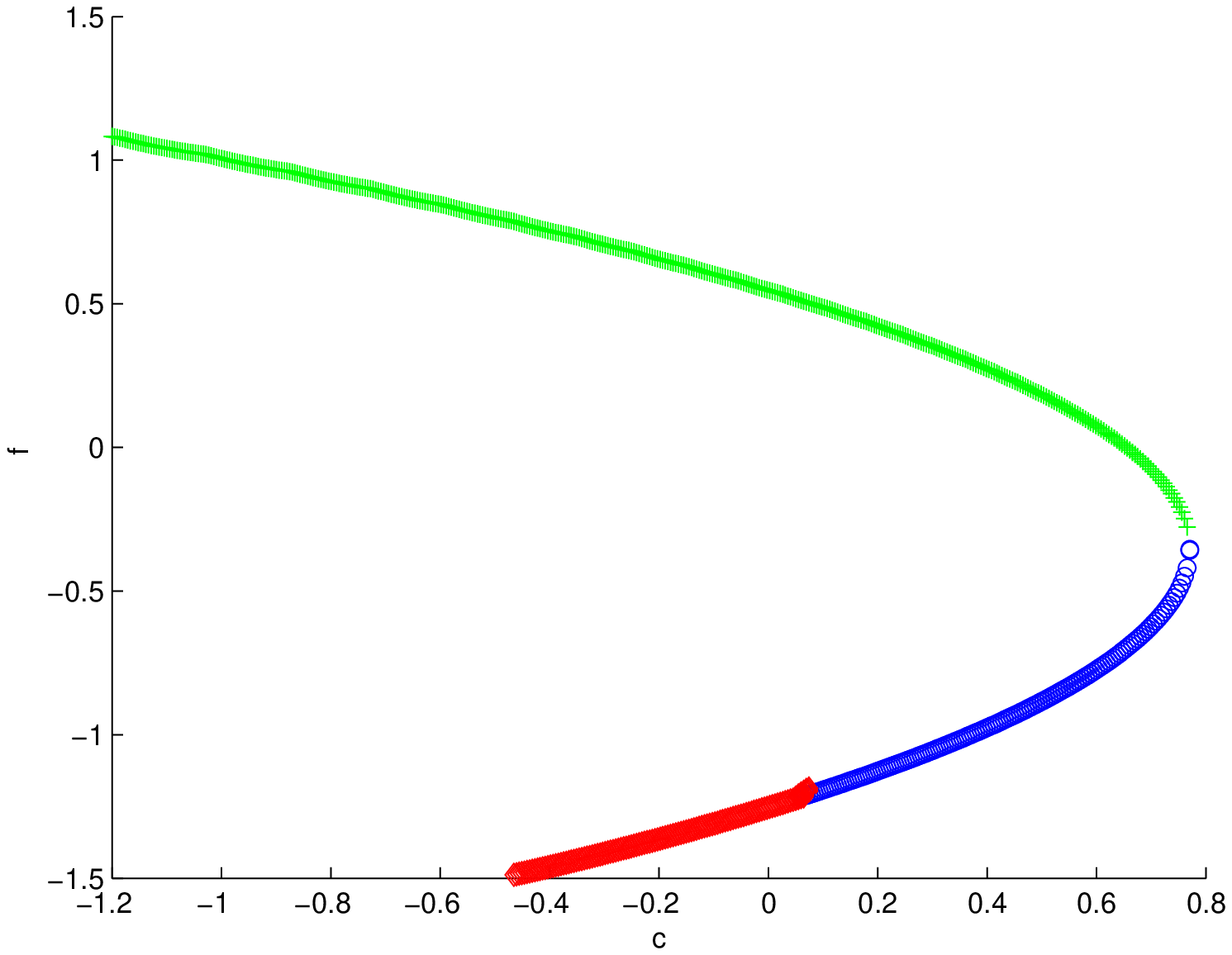}
\includegraphics[height=2in]{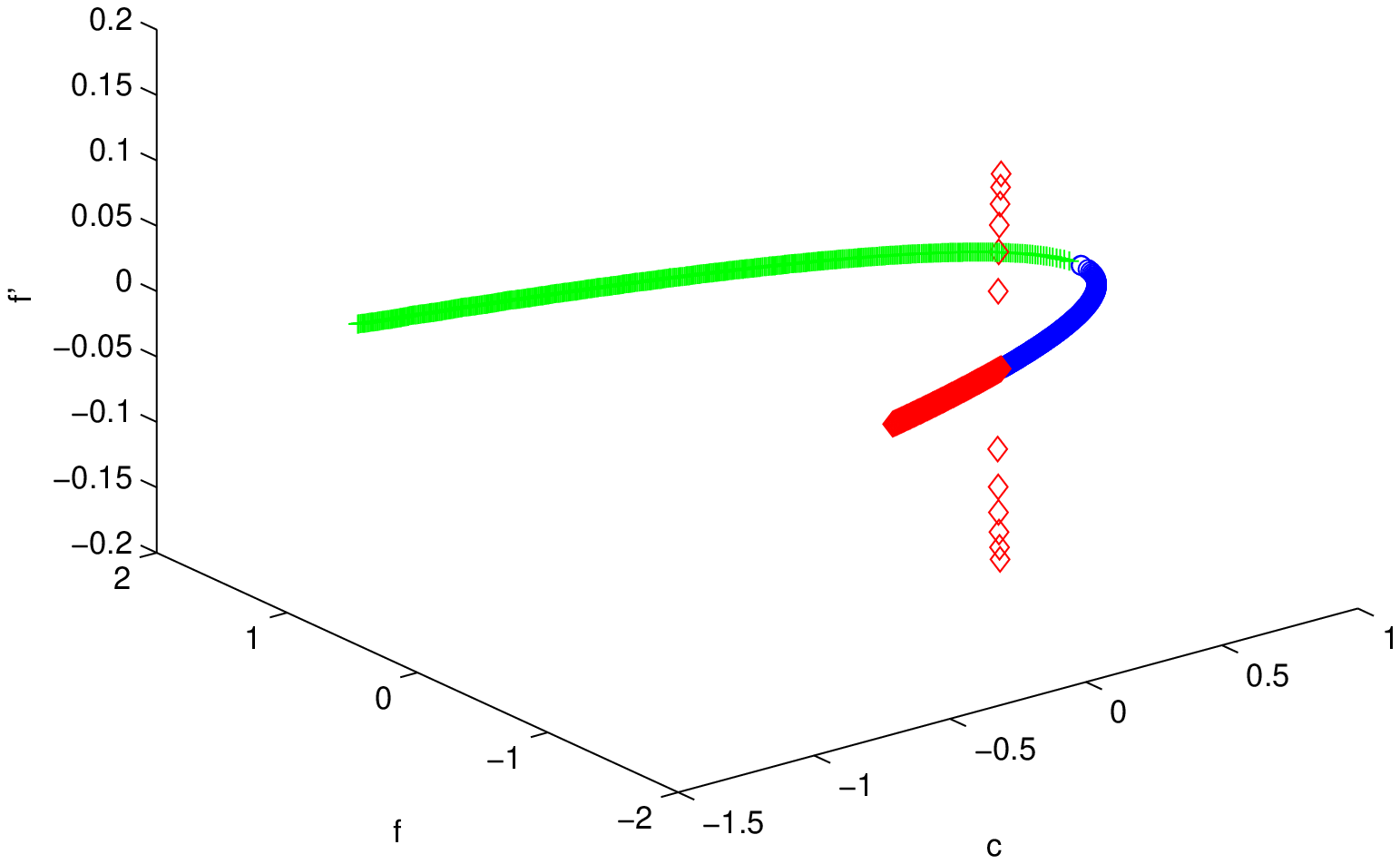}
\caption{Bifurcation diagram, coded by spectrum of
$\frac{d^2}{dx^2}-2f$: green = nonpositive spectrum, blue = one positive
eigenvalue, red = two positive eigenvalues}
\label{bif_diag}
\end{figure}

Consider the following spectial case of \eqref{limited_pde_unstable}

\begin{equation}
\label{example_pde}
\frac{\partial u}{\partial t}=\frac{\partial^2 u}{\partial x^2} - u^2
+ (x^2-c)e^{-x^2/2},
\end{equation}
where the choice of $\phi$ in \eqref{limited_pde_unstable} has been fixed.  The
bifurcation diagram for the equilibria of \eqref{example_pde} can be
found in Figure \ref{bif_diag}.  The bifurcation diagram is
parametrized by three variables: $c$, $f(0)$, $f'(0)$.  (Since the
equilibrium equation is a second-order ODE, it suffices to specify
each solution by its value and first derivative at 0.)  Based on the Theorem
\ref{unstable_thm}, the number of positive eigenvalues shown in Figure
\ref{bif_diag} corresponds exactly to the dimension of the unstable
manifold of each equilibrium.

\subsection{Frontier of the stable manifold}

\begin{figure}
\begin{center}
\includegraphics[height=1.33in]{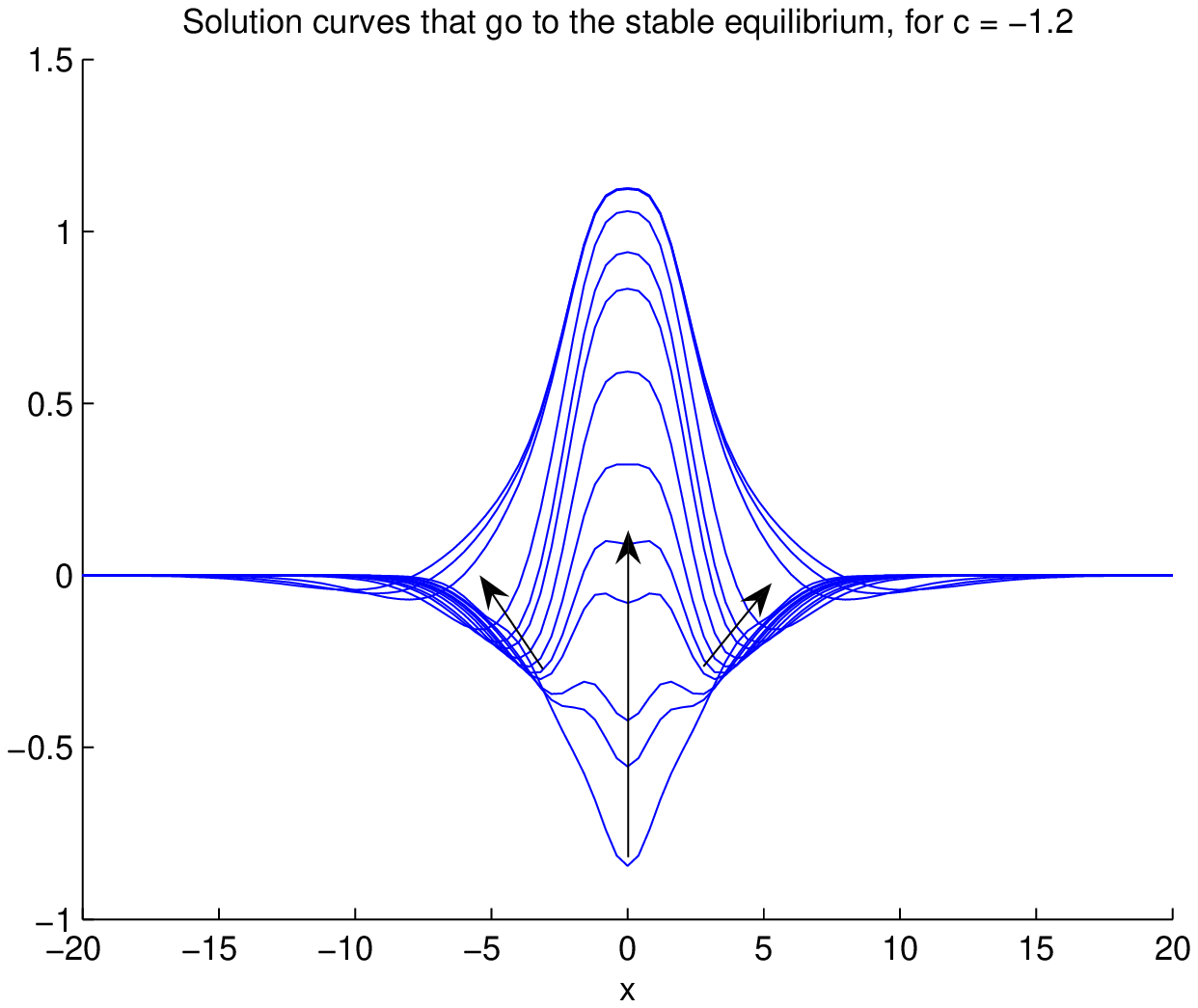}
\includegraphics[height=1.33in]{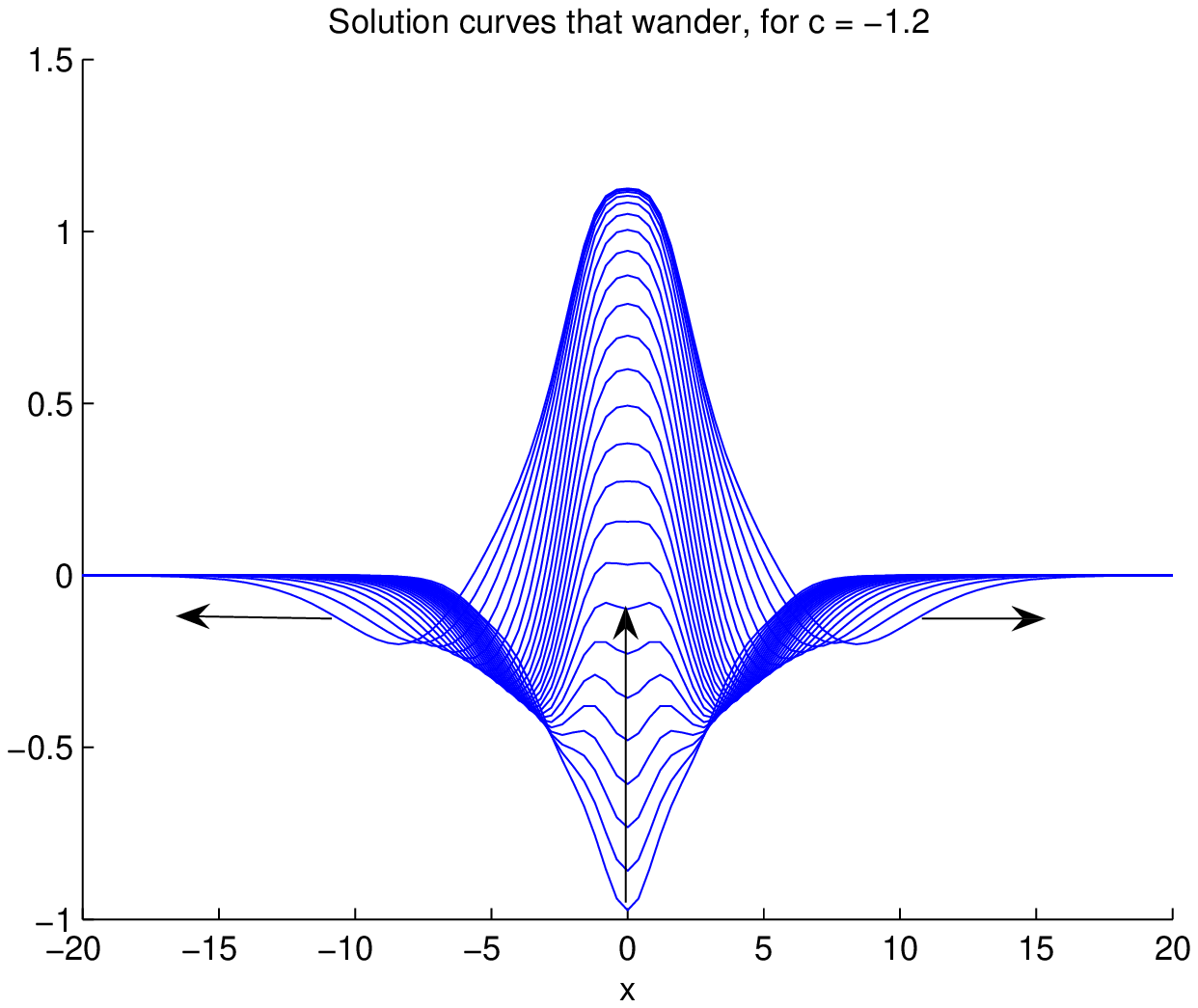}
\includegraphics[height=1.33in]{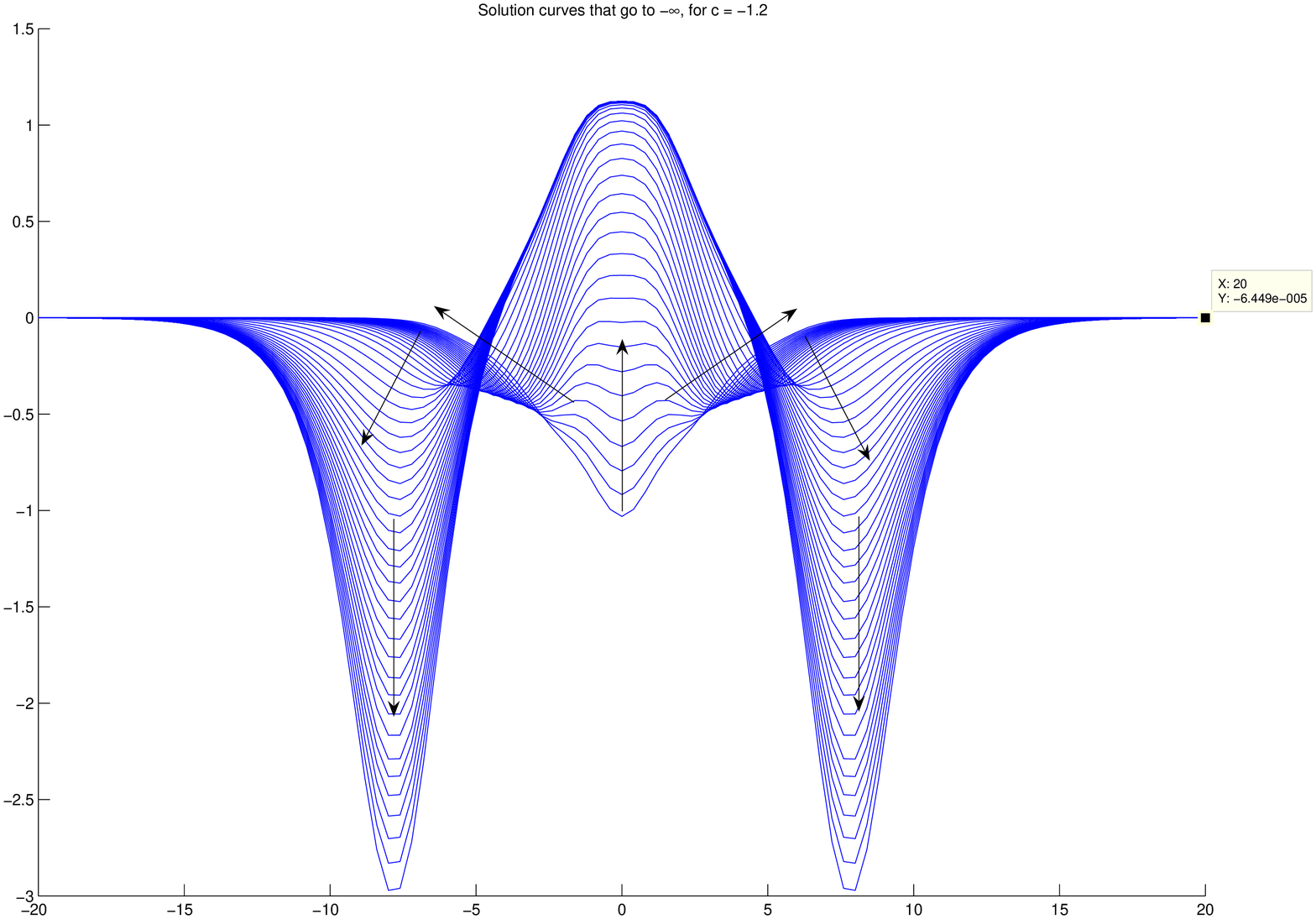}
\end{center}
\caption{Behavior of solutions near the frontier of the stable
  manifold of $f_0$ (horizontal axis is $x$)}
\label{frontier_fig}
\end{figure}

According to Figure \ref{bif_diag}, when $c=-1.2$, there is only one
equilibrium, $f_0$.  It has empty unstable manifold, though of course
it is asymptotically unstable (as is shown in
\cite{RobinsonInstability}).  On the other hand, $f_0$ has an infinite
dimensional stable manifold, which is not all of
$C^{0,\alpha}(\mathbb{R})$, as a consequence of the asymptotic
instability.  As a result, its stable manifold has a frontier in
$C^{0,\alpha}(\mathbb{R})$ (which may not be a boundary in the sense
of a manifold with boundary).  We are interested in the qualitative
behavior of solutions near and along this frontier.  We know by
Lemma 6 of \cite{RobinsonClassify} that if they tend to $f_0$
uniformly on compact subsets, then they do so uniformly.  It is
enlightening to use a numerical procedure to this end.  We start
solutions at the following family of initial conditions
\begin{equation}
\label{start_eq}
u_A(x) = f_0(x) + A e^{-x^2/10}.
\end{equation}
Using the Fujita technique (exactly as shown in 
\cite{RobinsonInstability}), we can show that for sufficiently negative $A$,
the solution started at $u_A$ will not be eternal.  As a result, the
family of initial conditions $u_A$ intersects the frontier of the
stable manifold of $f_0$.  An approximation to the value of $A$ which
corresponds to the frontier can be easily found using a binary search.
Some typical such solutions are shown in Figure \ref{frontier_fig},
and the approximate value of $A$ corresponding to the frontier is
$A\approx -2.15$

The qualitative behavior shown in Figure \ref{frontier_fig} indicates
that there is some kind of traveling disturbance in the frontier
solutions, which seems like a traveling wave.  However, such a
solution also appears to tend uniformly on compact subsets to $f_0$,
so in fact it converges uniformly.  (The uniform convergence is not
obvious from the figure, due to the numerical solution being truncated
at a finite time.)  The leading edge of this disturbance collapses to
$-\infty$ in finite time for solutions just outside the stable
manifold of $f_0$.

\subsection{Flow near equilibria with two-dimensional unstable manifolds}

\begin{figure}
\begin{center}
\includegraphics[height=2.5in]{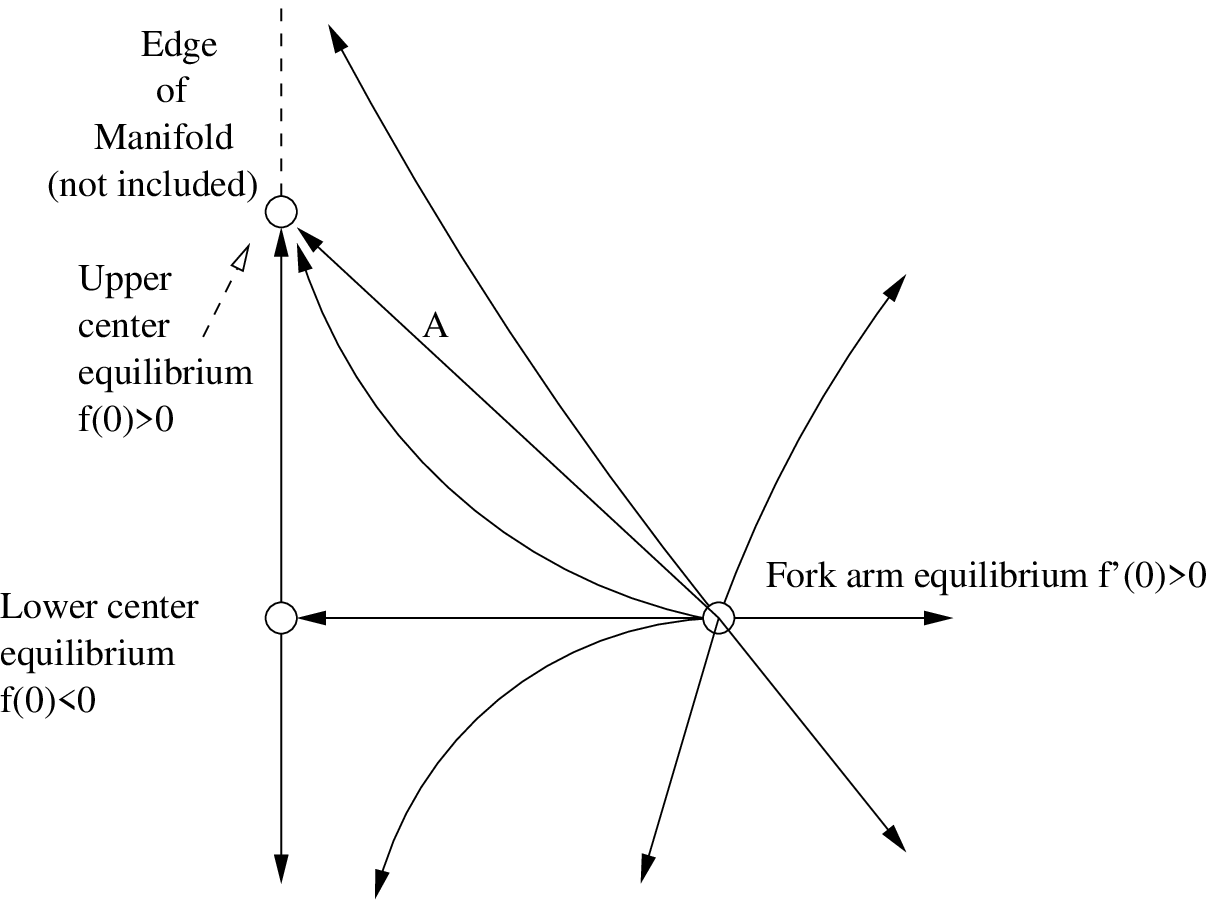}
\includegraphics[height=2.5in]{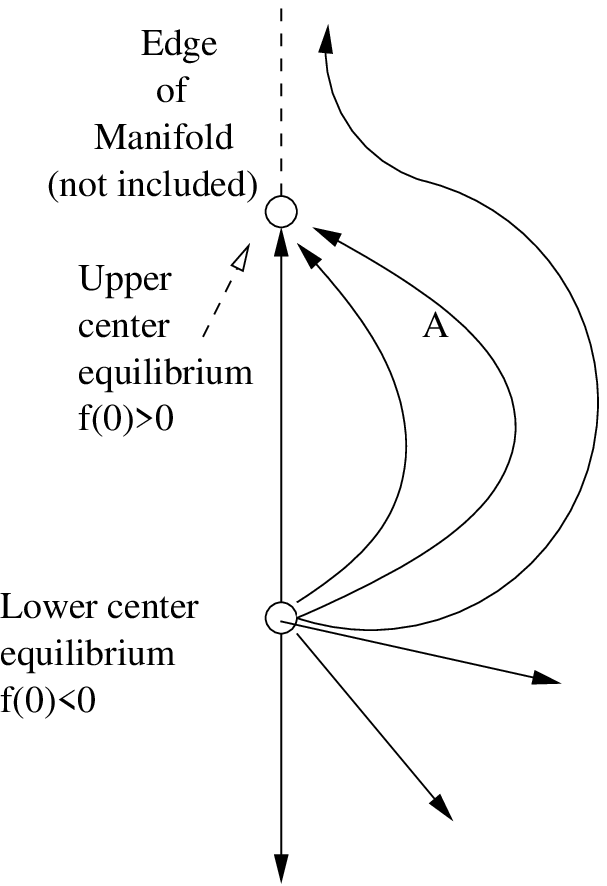}
\end{center}
\caption{Flow in the unstable manifold of a ``fork arm.''  $c=0.0600$
  (left); $c=0.0501$ (right)}
\label{flow_fig}
\end{figure}

\begin{figure}
\begin{center}
\includegraphics[height=2in]{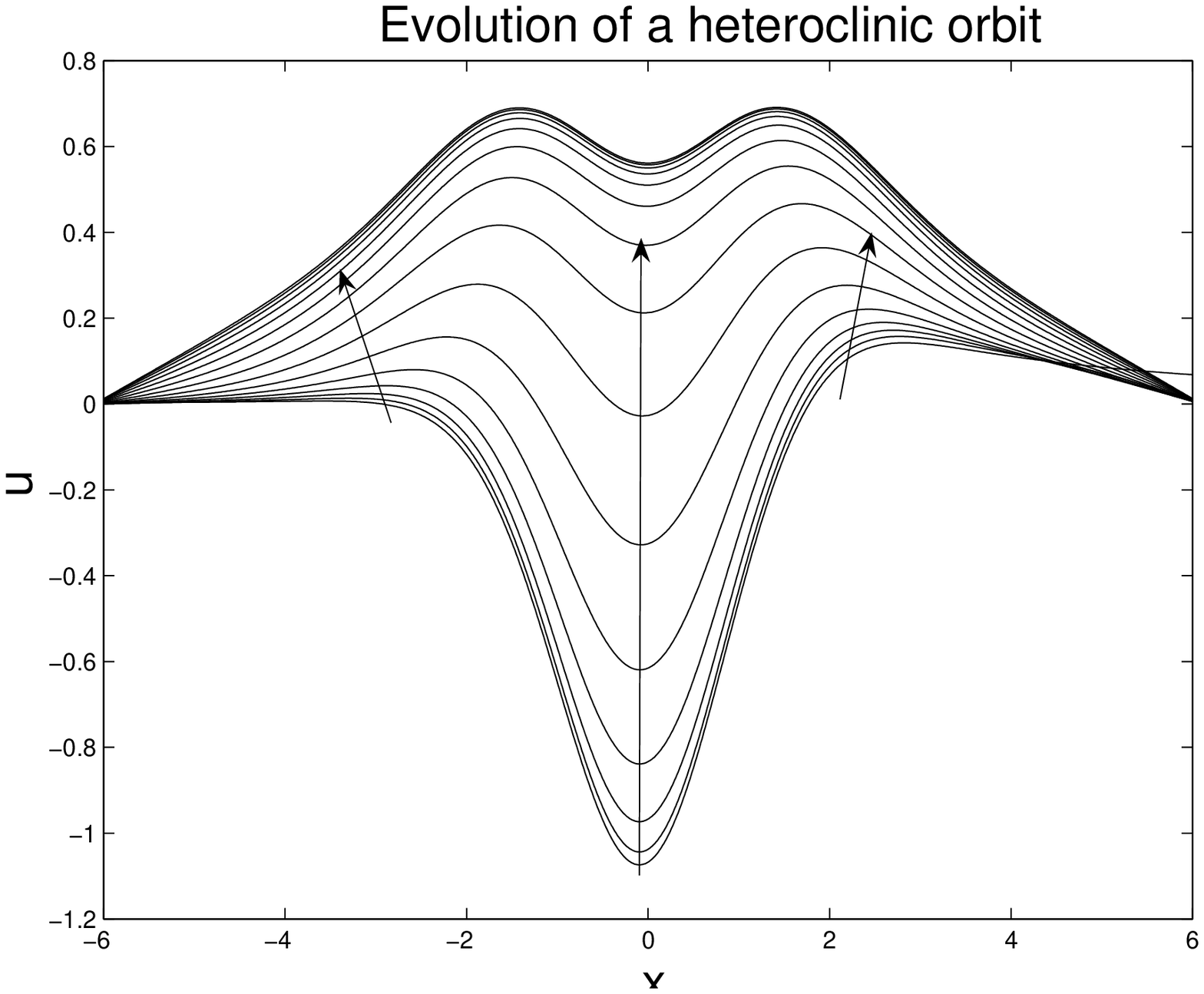}
\includegraphics[height=2in]{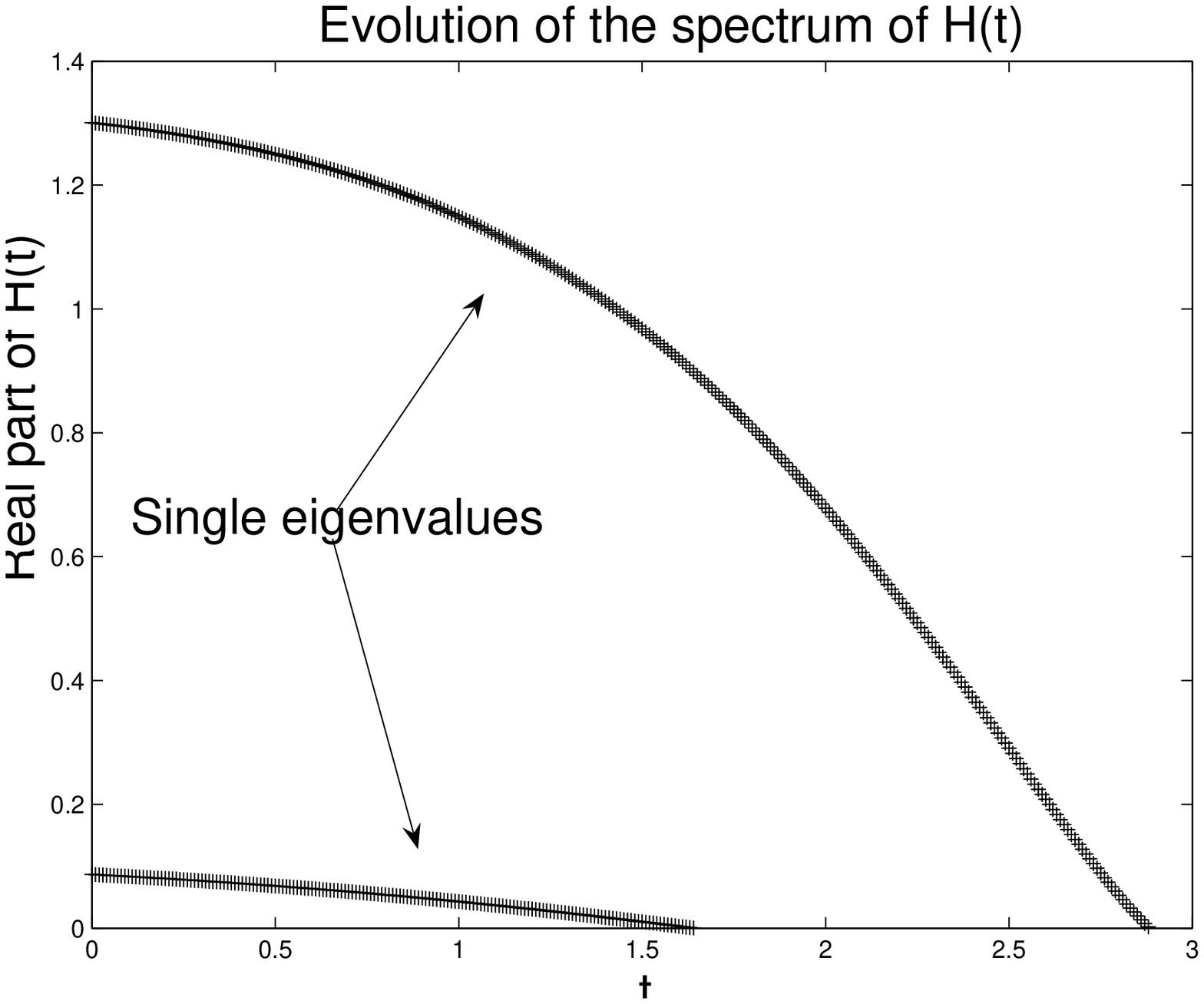}
\end{center}
\caption{A typical heteroclinic orbit to the left of boundary A, with
  the spectrum of $H(t)$ as a function of $t$.}
\label{solb_fig}
\end{figure}

Also of interest is the structure of the flow in the unstable
manifold of the ``fork arms'' which occur at $c=0.0740$, as they
approach the pitchfork bifurcation at $c=0.0501$.  Figure \ref{flow_fig} shows a
schematic of the flow based on numerical evidence.  Of particular
interest is the behavior near the boundary marked A.  Solutions to the
right of the boundary are not eternal solutions -- they fail to exist
for all $t$.  Solutions to the left of A are heteroclinic orbits
connecting the equilibrium with an unstable manifold of dimension 2 to the
equilibrium with an unstable manifold of dimension zero.  A typical
such solution is shown in Figure \ref{solb_fig}.

\begin{figure}
\begin{center}
\includegraphics[width=4in,height=2in]{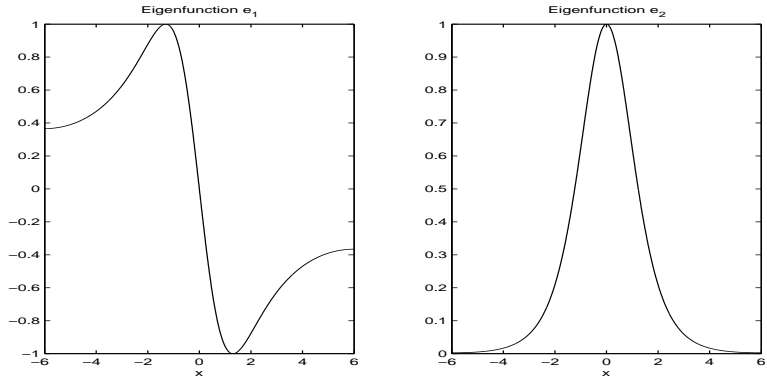}
\end{center}
\caption{Eigenfunctions describing unstable directions at $f_1$}
\label{eigenfuncs_fig}
\end{figure}

To examine solutions near the boundary A, we center our attention on
the case $c=0$, which has two equilibria, one of which (call it $f_1$)
has a 2-dimensional unstable manifold.  (This corresponds to the right
pane of Figure \ref{flow_fig}.)  If we linearize about $f_1$, the
operator $H=\frac{\partial^2}{\partial
x^2}-2f_1:C^{0,\alpha}(\mathbb{R})\to C^{0,\alpha}(\mathbb{R})$ has a
pair of simple eigenvalues, as is easily seen in the right pane of
Figure \ref{solb_fig} at $t=0$.  One of these eigenvalues is smaller,
to which is associated the eigenfunction $e_1$ in Figure
\ref{eigenfuncs_fig}.  The eigenfunction $e_2$ is associated to the
larger eigenvalue.  In Figure \ref{flow_fig}, $e_1$ corresponds to the
horizontal direction, and $e_2$ corresponds to the vertical direction.
From the proof of Lemma \ref{eq_findim}, it is clear that
$\{e_1,e_2\}$ spans the tangent space of the unstable manifold at
$f_1$.  Therefore, we specify initial conditions $u_{A,\theta}(x)$ for
a numerical solver using
\begin{equation}
\label{ic_eq}
u_{A,\theta}(x)=f_1(x) + A \left(e_1(x)\cos \theta + e_2(x) \sin\theta\right).
\end{equation}

\begin{figure}
\begin{center}
\includegraphics[width=5.75in]{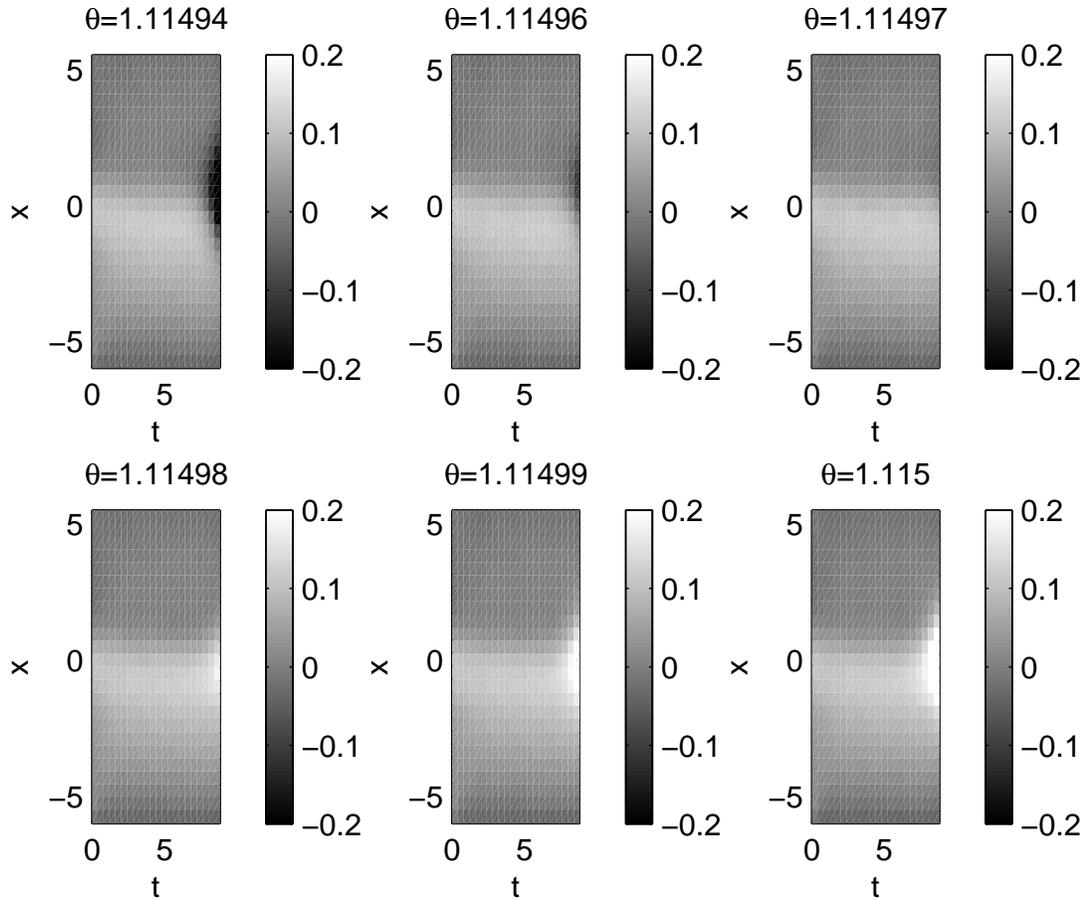}
\end{center}
\caption{Difference between equilibrium $f_1$ and the numerical
  solution started at $u_{A,\theta}$, where black indicates a value of
  -0.2, and white indicates 0.2.  The horizontal axis represents $t$,
  and the vertical axis represents $x$.  $A=0.1$ in all figures.
  Starting from the upper left, $\theta=1.11494, 1.11496, 1.11497,
  1.11498, 1.11499, 1.115.$ }
\label{perturb_fig}
\end{figure}

(Taking $A$ small allows us to approximate solutions which tend to
$f_1$ in backwards time.)  Since the perturbations along $e_1,e_2$ are
quite small, and indeed the eigenvalue associated to $e_1$ is much
smaller than that associated to $e_2$, examining the numerical results
of evolving $u_{A,\theta}$ is quite difficult.  The behavior along the
boundary occurs at a much smaller scale than $f_1$, yet is crucial in
determining the long-time behavior of the solution.  To remedy this,
the boundary behavior is better emphasized by plotting
$u_{A,\theta}(t,x) - f_1(x)$ instead.  Figure \ref{perturb_fig} shows
the results of evolving initial conditions \eqref{ic_eq} for $A=0.1$
and various values of $\theta$.

Solutions in Figure \ref{perturb_fig} show a similar kind of behavior
as in the case of the frontier of $f_0$.  There is a traveling front,
which moves very slowly in the negative $x$-direction.  However, the
behavior is quite a bit more delicate.  The determining factor in
locating the frontier of $f_0$ is the perturbation in a direction
roughly like $e_2$, which has a large eigenvalue.  On the other hand,
for $f_1$, Figure \ref{flow_fig} indicates that such a direction is
not parallel to the boundary of the connecting manifold.  (The
boundary direction is some linear combination of $e_1$ and $e_2$, with
a numerical value for the angle $\theta$ being roughly 1.114975
radians.) The eigenvalue associated to $e_1$ is roughly ten times
smaller, and therefore perturbations in that direction are much more
sensitive.  Additionally, the action of the flow is therefore
primarily in the direction of $e_1$, which tends to mask effects in
other directions.  For this reason, it was visually necessary to
postprocess the numerical solutions by subtracting $f_1$ from them.
Otherwise the presence of the traveling front was unclear.

\section{Conclusions}

We have shown that the tangent space at an equilibrium splits into a
finite dimensional unstable subspace, and infinite dimensional center
and stable subspaces.  However, it is quite clear by 
\cite{RobinsonInstability} that the center subspace is nonempty and large.
Indeed, considering the work of \cite{Souplet_2002}, the center and
stable subspaces are not closed complements of each other.
Additionally, we have given conditions for the space of heteroclinic
orbits to have a finite dimensional cell complex structure.

\bibliography{cellcomplex_bib}
\bibliographystyle{plain}

\end{document}